\documentclass[12pt]{amsart}
\usepackage{amsfonts}
\usepackage{graphicx}
\usepackage{times}
\usepackage{tipa}
\newcommand{\dis}{\displaystyle}

\newcommand{\Nn}{\mathbb{N}}
\newcommand{\Zz}{\mathbb{Z}}

\newcommand{\beq}{\begin{eqnarray*}}
\newcommand{\eeq}{\end{eqnarray*}}
\newcommand{\ord}{\mbox{\rm ord}\, }
\newcommand{\arf}{\mbox{}^\ast \! }
\newcommand{\pp}{\mathfrak{p}}
\newcommand{\qq}{\mathfrak{q}}

\newcommand{\al}{\mbox{\textscripta}}
\begin{document}
\title{A Scientific Biography of Cah{\.i}t Arf (1910-1997)}
\author{Al{\.i} S{\.i}nan Sert\"{o}z}
\address{Bilkent University, Department of Mathematics, 06800 Ankara, Turkey \\ {\tt http://www.bilkent.edu.tr/$\tilde{~}$sertoz}}
\email{sertoz@bilkent.edu.tr}
\date{16 January 2013, \textbf{version 1.22}}

\begin{abstract}
In this paper we analyze both the scientific activities of Cahit Arf, a Turkish mathematician, and the social context in which he worked. We also discuss his work and social environment leading to the discovery of Arf invariant, Arf rings, Arf closure and Hasse-Arf theorem.
\end{abstract}

\maketitle
\thispagestyle{empty}

\newpage
\thispagestyle{empty}
\mbox{} \\
{\small \it this page is left blank}
\newpage
\setcounter{page}{1}
\thispagestyle{plain}

\section{Introduction}
Cahit Arf (1910-1997) was a  Turkish mathematician with a legendary fame in his country. At the international mathematical circles he was known via certain concepts bearing his name. The most well known of these is probably the  \emph{Arf invariant} of algebraic topology. It even became a subject of humor as depicted in a cartoon in a proceedings book\footnote{Jig/reel image, taken from page xiii of Current trends in algebraic topology. Part 1 (Proceedings of the Conference held at the University of Western Ontario, London, Ont., June 29--July 10, 1981, edited by Richard M. Kane, Stanley O. Kochman, Paul S. Selick and Victor P. Snaith)} as follows:
\begin{center}
\resizebox{11cm}{!}{\includegraphics{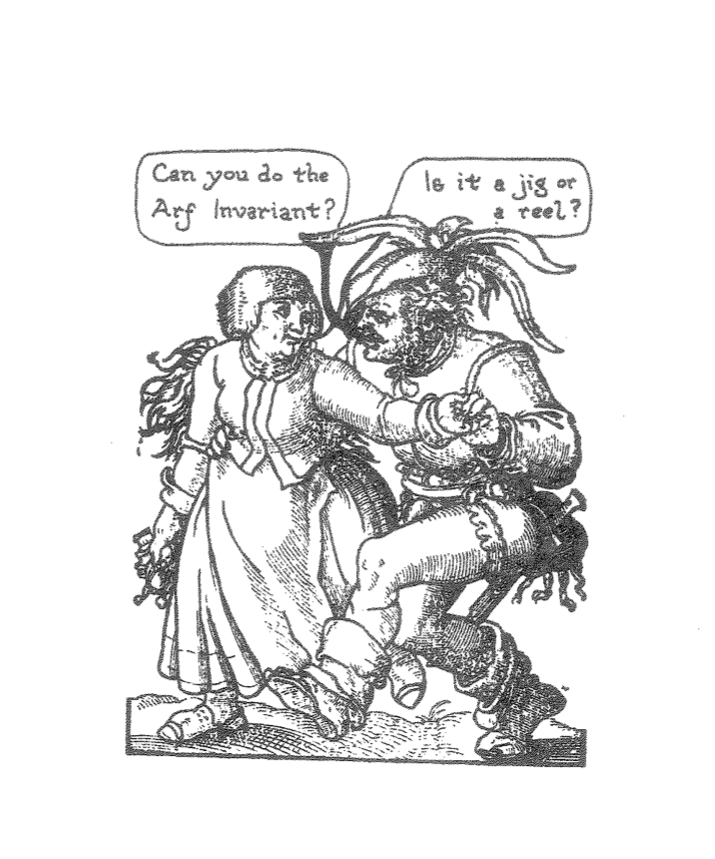}}
\end{center}

Arf invariant also appears in the Kervaire-Arf invariant problem. For an even dimensional manifold, the cup product on the middle cohomology group with $\Zz/2\Zz$ coefficients defines a non-degenerate bilinear pairing. The Arf invariant of this pairing, which is either one or zero, is known as the Kervaire invariant of the manifold. The Kervaire-Arf invariant problem is to determine if there are manifolds with Kervaire invariant one. Throughout the long history of this problem it was known that if the manifold has Kervaire invariant one, then the dimension of the manifold has to be of the form $2^{2k+1}-2$ for $k>0$. It was also known that such manifolds do exist for $k=1,\dots,5$. Recently M. A. Hill, M. J. Hopkins and D. C. Ravenel proved that no such manifolds exist for $k\geq 7$, see \cite{hill1,hill2,hill3}. The case $k=6$ remains open.

In Turkey the formula of a special case of Arf invariant, together with Arf's portrait, appears on the back of a 10 Lira banknote. This note is tender since 2009.
\begin{center}
\resizebox{11cm}{!}{\includegraphics{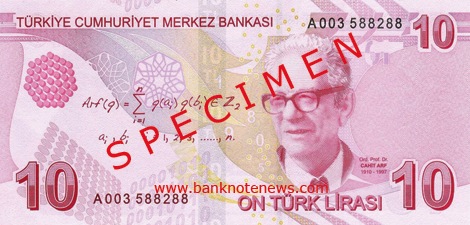}}
\end{center}
Moreover in 2001 a commemorative coin\footnote{Two series of commemorative coins were issued, one in 2001 and the other in 2003. The denomination of Arf coin is 7500000TL with 2001 currency. It was then worth slightly above 6USD.} was issued in his name as a collectors item.
\begin{center}
\resizebox{5cm}{!}{\includegraphics{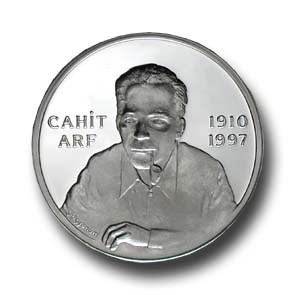}}
\resizebox{5cm}{!}{\includegraphics{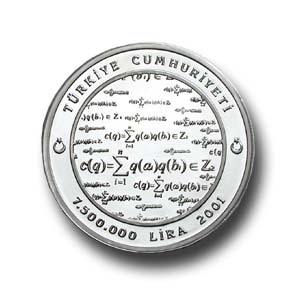}}
\end{center}

Arf's life in Turkey spanned periods of great change and numerous social turmoil. As most intellectuals witnessing such changes in their societies he developed strong views on these issues and did not refrain from being outspoken about his ideas. The result of this was that he became a well known figure in public life.  Through his public fame people became aware of mathematics.

The purpose of this article, written so many years after his death, is to give an analysis of Cahit
Arf's scientific production, to describe the social
and humane context in which this work progressed and discuss his influence on his contemporaries and students.

\section{Life and Influence of Arf}
Cahit Arf was born on 11 October, 1910 in Selanik (now Thessaloniki, Greece) which was then an Ottoman town on the Balkan peninsula. With the outbreak of the Balkan war in 1912 Arf family moved to Istanbul\footnote{It must be noted that the Ottoman tradition did not use family names. People were identified by their father's name together with their own name.  It was only in 1934 during the Turkish Republic era that a law was passed to bring order to this situation and surnames were institutionalized. Arf's family then decided to create and use the word \emph{Arf}, obtained  from \emph{Arif}=wise. This was to avoid inadvertently obtaining a name also chosen by other families. They succeeded in choosing a unique name in Turkey but later there were numerous inquiries from Arabic countries asking if they had any relations with a certain Arf Pasha.}.

Cahit Arf's childhood encompassed the Balkan wars, the  World War I, the grand war at Gallipoli, the Greek invasion of western Anatolia and the invasion of Istanbul by the Allied Powers. Ottoman Empire was dismantled with the S\`{e}vres Treaty\footnote{This treaty of 1920 was never adopted. It was superseded by the Lausannes Treaty in 1923.}. Mustafa Kemal and friends rebeled against the occupying powers in 1919. Cahit Arf was 9 years old when his father, an employee in the postal office, decided to move to Ankara and join the ranks of this movement.

When finally Turkey emerged as a new independent parliamentary republic in 1923 Cahit Arf was 13 years old.  The new Republic was poor and deplete of any material sources due to countless wars fought. However the new republic managed to induce hope in the people. The victory of the last battle fought and won was enough cause to be full of hope and feel invincible.  There was one way to go and that was to believe that everything could be rebuilt from scratch.  This was the environment while Arf was becoming a teenage. This attitude affected him so much that in his mathematical articles, to be written years later, one would witness how he manoeuvres the treatment of the problem to a phase where he has to develop all the fundamental tools he needs for the final blow, as opposed to the more common manner of reducing the parts of the problem to already known cases.

His mathematical talent gave its first signals in public school when young Cahit was excelling in grammar. His success at grammar shaded his ease in mathematics but soon a young teacher discovered his natural aptitude towards mathematics.  With the guidance of this teacher Arf  studied  through the \emph{Elements} of Euclid, \cite{CA}.

After primary school, an overseas economic crises helped to shape Arf's high school education.
In 1926 devaluation of French Frank against Turkish Lira provided his father with the financial means of sending him to France. It was a popular place for the Young Turks of the late Ottoman era and besides France was the first Allied Power to sign a treaty with the young Republic. Cahit Arf was thus sent to the prestigious St Louis Lyc\'{e}e. His weakness in the French language was more than made up with his extraordinary grades in mathematics exams and he graduated in two years instead of the expected three years. This was a school well known for its intention to prepare students for \'{E}cole Normale Sup\'{e}rieure and Polytechnique.

Upon returning to Turkey Arf obtained a state scholarship to continue his studies in \'{E}cole Normale Sup\'{e}rieure  in France. His stay at \'{E}cole Normale  coincides with that of Jean Leray who was there between 1926 and 1933. We do not know if they met but it is safe to presume that they met and made friends\footnote{Leray's name was mentioned to Murad {\"{O}}zayd{\i}n after a complex analysis final exam.
Arf asked the proof of Leray's Theorem  on page 189 of Gunning and Rossi's book \cite{gunning}.
When confronted, by {\"{O}}zayd{\i}n who wrote the exam, with the
inappropriateness of such a question in undergraduate level, he said he asked this question because he knew
Leray from his France years and he was curious about what his friend did.}.  He must be the hero in the following story. Arf was in the habit of encouraging young mathematicians to discover mathematics by themselves rather than learn it from written sources. This is certainly a good advice provided that it is taken with a grain of salt, but Arf's insistence on it was at times beyond reason. To support his cause he would tell, with some exaggeration,  how in his \'{E}cole Normale years in Paris, he would never attend classes but would instead ask his friend what the class matter was that day and then proceed to develop that theory himself. It is very unlikely that his friend was an ordinary student. Leray, our candidate for this position,  at that time was a senior student and was already in command of the lecture matters. It is probable that Leray was auditing a certain interesting but introductory course which Arf was officially taking as a junior student and it was only with that course and with none other than Leray as his friend that Arf accomplished this feat of which he was so proud later on. With a lesser friend the project could not have taken off the ground. Who else can summarize the course material so clearly and correctly for a young budding student to take it from there and work his way towards a satisfying conclusion!

Nonetheless this attitude of developing everything from basic principles stayed with Arf throughout his mathematical life.

Arf stayed at \'{E}cole Normale Superior between the years 1928  and 1932 and returned to Turkey
as the first Turk to graduate from this prestigious
school. This was later acknowledged  in 1994 when a medal\footnote{Commandeur dans l'Ordre des Palmes Academiques.}  from France was presented to him during a ceremony at Palais de France, the French consulate in {\.I}stanbul.

In 1932 as a young  graduate  Arf decided to dedicate his life  to teaching mathematics at high schools. He was thinking of going to Kastamonu\footnote{Today Kastamonu does not enjoy the same glory it once had during the early years of the Republic. In particular Kastamonu Lyc\'{e}e then was very well known and reputed for the quality of its graduates. It was probably this school that was in Arf's mind. He must have been influenced by the grandeur of the city and the school at that time when he had a forced stop over there on his family's migration to Ankara in 1919 to follow Mustafa Kemal.} but was convinced to stay at {\.I}stanbul and lecture at Galatasaray Lyc\'{e}e\footnote{This is a school whose origins date back to 1481. It has undergone several revisions and today functions as a prestigious high school with French as the medium of education.}. At that time high school education was extremely strong in Turkey and Arf's decision to teach in a high school, given his enthusiasm and idealism, is understandable. In fact Oktay Sinano\u{g}lu, the youngest scientist to ever receive professorship at Yale University\footnote{He was 26 when he was promoted as a professor at Yale University in 1962. Since then he has been several times a Nobel nominee  in chemistry.}, repeatedly claims that his initial success at Yale was due to his high school education in Turkey. University was not yet a popular concept in Turkey then. However the new Republic was set to follow the modern world. Reforms were made to universities in 1933 and Arf was called on to duty as a young instructor at {\.I}stanbul University\footnote{His rank was \emph{candidate for docent}, probably the equivalent of assistant professor of today.}.

Incidentally there was another prominent Turkish mathematician, Kerim Erim\footnote{He is the maternal grandfather of G\"{u}lsin Onay, international concert pianist and State Artist of Turkey.}, on the board who planned and executed these reforms. Kerim Erim is credited to be the first Turkish mathematician with a PhD degree. He received his degree in 1919 from the Frederich-Alexanders University in Erlangen, Germany.\footnote{See \cite{alp}, page 437.} He was instrumental in attracting Cahit Arf, Ferruh \c{S}emin and Ratip Berker to university life. All these gentlemen would later become leaders of the Turkish mathematical community and guide the new generations towards excellence in academic world.

It did not take Cahit Arf long to realize  that he needed to do  graduate study  in mathematics. This is how he arrived in 1937
at G\"{o}ttingen to study with Helmut Hasse. Arf was thinking about classifying all algebraic equations which could be solved by radicals. This would involve Galois theory and the classification of all solvable groups already done by Jordan. The sheer enormity of this problem is a reflection of that self confidence, idealism and determination of the young Turkish Republic of the time, which Cahit Arf inadvertently radiated himself. Hasse naturally guided his young student towards a realistic subproblem. Years later in 1974 during a symposium in Silivri, Turkey, Hasse would recall that after taking his dissertation problem Arf had disappeared from the scene for a few months only to come back with the solution.  After developing this solution with his advisor, Arf published his findings in Crelle's Journal, \cite{carf1}. This theorem is now known as Hasse-Arf theorem.

It would be interesting to know what Hasse would write in his \emph{Mathematischen Tageb\"{u}cher} about his mathematical encounter with Arf. Unfortunately Hasse kept his mathematical diary only between 1923 and 1935.

Those were the  days before the World War II in Europe. There was instability and turmoil in the western world. However Hasse was concerned that times in Turkey would be even more unbearable. He asked his student to stay for another year and suggested he work on a new problem, classification of quadratic forms over a field of characteristic two. This was a problem which Witt had solved recently for other characteristics. Arf solved the problem for characteristic two and published it in Crelle's journal, \cite{carf2}. In that article he totally classified quadratic forms over a field of characteristic two. His findings were  algebraic in nature but they were so fundamental that they later found wide applications in the classification problems of topology. His findings there are known today as \emph{Arf invariants}. As we mentioned above, this is also the source for the Kervaire-Arf invariant problem.

During the Second World War, Turkey  benefited from the visits of several  German scientists who escaped discrimination. For example mathematician Richard von Mises was in Turkey between 1933 and 1939. In fact it was through the recommendation of von Mises that Arf went to study with Hasse, see \cite{alp}. Another German mathematician visiting Turkey at this time was  William Prager who arrived in Turkey in 1934 and stayed until 1941 when he immigrated to USA. For a detailed narration of German mathematicians in exile in Turkey during this era and their influences on Turkish mathematicians see the recent work of Eden and Irzik \cite{alp}.

Meanwhile the mention of Hasse's alleged intimacy with the Nazi authorities later always hurt Arf who knew Hasse so closely. His only defense would be to remark that ``Hasse was a very kind man". In fact the friendship between these two mathematicians grew over the years to encompass their families. Today Peter Roquette has recovered letters written in Turkish between Mrs Arf and Hasse  for the purpose of giving Hasse an opportunity to practice his Turkish. Hasse has visited {\.I}stanbul several times on Arf's invitation. Tosun Terzio\u{g}lu has recently told me that he learned from Mrs Arf that on one of these visits Arf came home with a broad smile on his face. He explained to his wife that Hasse told him that day  that he can now address him with the more friendly pronoun \emph{du}, the second personal singular pronoun in German.

An English mathematician, with the name of Patrick du Val\footnote{Du Val was then a bright mathematician. He was sent to {\.I}stanbul by the British Council probably to balance the effect of German scientist on the young Republic. He stayed at {\.I}stanbul University until 1949, learned Turkish and delivered his lectures in Turkish. He even wrote  mathematical textbooks in Turkish. He later became famous for his work on singularities, some of which was carried during his residence at {\.I}stanbul. In fact today certain types of singularities bear his name in algebraic geometry. Du Val returned to {\.I}stanbul University in 1970 after his retirement and stayed there for another three years.}, also visited Turkey during that period. In one seminar he talked about the problem of determining the multiplicity sequence of a space curve branch, starting with the initial data. The problem was then totally understood for plane curves. And for space curve branches Du Val himself had found some geometrically significant infinitely near points which were sufficient to determine the multiplicity sequence. The problem was  to describe Du Val's significant points starting from the formal parametrization of the branch.  Arf totally solved this problem and the rings he constructed in his solution and the process of obtaining those rings are now known as \emph{Arf rings} and \emph{Arf closure}.

The algebraic structures associated to geometrically singular points in algebraic geometry generally contain gaps compared to the corresponding structures associated to smooth points. Resolving the singularity generally involves filling these gaps algebraically and hoping that the geometric counterpart of the process eventually smoothes out the given singularity. This phenomenon is also used in Hironaka's famous work on resolution of singularities \cite{hironaka}, and this is exactly what Arf did to determine the multiplicity sequence which measures the \emph{badness} of the singularity at each step of the resolution process.

This solution has an interesting story which Cahit Arf himself narrated in detail and in full sincerity to me in Gebze, at the Fundamental Sciences Institute,  around 1985. When Du Val was talking about the resolution of branch singularities in the seminar, back in 1944, Arf was listening but had no idea about what a blow-up was. However he was observing certain algebraic regularities in the so called blow-up process. To cover up his ignorance he claimed  that Du Val was too much involved in the geometry of the problem. He even claimed that the problem was probably totally algebraic. Upon Du Val's challenging him to solve the problem without geometry, he stayed home for a week, due to a severe cold, and returned to the department with the solution, \cite{carf8}. That paper is a force de tour of developing all the tools \emph{by hand}, which is typical of an Arf manuscript.

Meanwhile his successful career at {\.I}stanbul University brought him a corresponding membership of the Academy of Sciences in Mainz\footnote{Search Google for ``Kategorie:Mitglied der Akademie der Wissenschaften und der Literatur".} in Germany. He went to Mainz
to deliver a lecture on the arithmetic analogue of
the Riemann-Roch theorem. Afterwards he stayed a month in Hamburg on the invitation of Hasse. There he gave a series of lectures on Riemann-Roch theorem on number fields. His paper on this topic appeared in the proceedings of Mainz Academy, \cite{carf20}. Peter Roquette who was then in Hamburg attended Arf's lectures and comments\footnote{Private communications.} that this particular piece of work deserves to be mentioned among Arf's mathematical achievements.

During this time in {\.I}stanbul, Arf started to talk with engineers and physicists. One of them was Mustafa {\.I}nan who was one of the first scholars sent abroad for graduate work\footnote{He had his Ph.D. at ETH Z\"{u}rich in 1941.}. He did his graduate studies in Switzerland on the forces acting on a collapsing bridge. He found some stress lines which affected the bridge's tolerance. His work was mostly experimental. Later this became an ideal problem for Arf to stretch his muscles on. He set out to give a mathematical modeling of these stress lines. In a series of papers he managed to solve the problem in two dimensions, \cite{carf9, carf10, carf12, carf13, carf16, carf17}.  The beauty and the strength of his work received for him in 1948 the first and for that matter the only {\.I}n\"{o}n\"{u} award given.
In fact he received the award for his first article.
The prize money was so generously set that he could pay off his mortgage on his house which was on the banks of Bosphorus where he lived for the rest of his life. He never returned to the elasticity problem again after this series of articles.

Encouraged by his success in applied mathematics he later attacked some problem on statistical mechanics. He was not as successful as on the elasticity problem but he enjoyed his collaboration there with the young scientists of the day, \cite{carf21}.

In 1949-1950 Arf visited Maryland University. Here history repeated itself. Recall that the work leading to the discovery of Arf rings had started in a seminar which Du Val was giving. Arf had then insisted that the problem would accept a totally algebraic modeling and their discussions with Du Val resulted in Arf's bringing out the structure he claimed to have existed.
He would openly admit himself that he was not a good reader of mathematical texts and would cherish at seminar discussions. Another such discussion which led to a theorem took place as he was in Maryland. He would later narrate that he was attending this seminar on a certain operator about which Arf knew practically nothing. But he observed that there was an algebraic pattern beyond these discussions and through his interactions with the people in the seminar he managed to give an algebraic presentation of those operators, see \cite{carf11}. Yet another such work which had its origins at Maryland is \cite{carf15}.

He was not isolated from the society as one would anticipate in the evidence of his mathematical achievements. He was actively interested in social matters and he was a well known figure in society. In fact he was sometimes invited to talk at meetings which are not necessarily scientific in the sense he was used to. We see an example of this in the speech he gave in Hatay on possible applications of mathematics in economical affairs, back in a 1956 meeting, \cite{arf-hatay}\footnote{I am grateful to Ali Erhan \"{O}zl\"{u}k for sending me a copy of this article.}

After retiring from Istanbul University in 1962 he worked for a year at Bo\u{g}azi\c{c}i University which was then called Robert College, in Istanbul. He was invited to the Institute for Advanced Studies as a visitor. He wrote to the institute to ask what was expected of him during his visit. He would later tell his students that the institute asked nothing of him besides his being there and pursuing his own interests. This was in accord with his conviction about how science should be done and he could now show a strong ally in the name of the Institute. He visited Princeton between 1964-1966.

After Princeton he visited California, Berkeley for a year and was about to accept an offer from H.S.M. Coxeter to visit Toronto when he received a letter of invitation from Erdal {\.I}n\"{o}n\"{u}\footnote{Erdal {\.I}n\"{o}n\"{u} is a physicist and received the Wigner Medal in 2004. He is the son of the late president {\.I}smet {\.I}n\"{o}n\"{u} for whose name the {\.I}n\"{o}n\"{u} award was given to Arf back in 1948. {\.I}smet {\.I}n\"{o}n\"{u} was a friend of Atat\"{u}rk and the second president of the Republic.},  who was writing on behalf of the  president of METU, the Middle East Technical University. METU was founded in 1956 and was aiming to reach world standards and was therefore keenly interested in recruiting prominent Turkish scientists. In his letter Erdal {\.I}n\"{o}n\"{u} not just included a polite invitation but also attached a plane ticket to Ankara. This was the beginning of the golden age of METU mathematics department.

Arf's residence at METU mathematics department did not produce any significant mathematical output on his part. But he was an idol and the aura he created was the cradle in which present day Turkish mathematics flourished. The writer of these lines still remembers with nostalgia the days when Arf was teaching undergraduate algebra at the blackboard and they were calculating first cohomology groups of group extensions.
Those however were ironically the most difficult times for Turkey. METU suffered several blows from the political powers due to student riots.  At the most crucial times Arf stood with determination for METU faculty and students risking his reputation and comfort.

A certain dialogue with the governing powers of the country at that time is still remembered today with awe and admiration. Following a criticism from government authorities about the education at METU,  he took word and boldly expressed that university is not the place where you teach what is known but it is the place where you teach how to reach what is not known. He followed his point to the end and stayed critical of all reforms to bring standards to university education, see for example \cite{arf-yok}.

He was also instrumental in the establishment of T\"{U}B{\.I}TAK in 1963, the Turkish Scientific and Technical Research Council. He served as the head of the scientific committee for some time. But Arf's relaxed approach to matters of the world was not the most efficient way to run business in most cases. He was a bad administrator, one claim he would not bother to disagree. In fact he advocated that good scientists should stay away from office.
He once even put this idea to practice at METU mathematics department, with disastrous results.
His students following his advice on mathematical matters did not follow this dubious wisdom and they did well later as head administrators of T\"{U}B{\.I}TAK and also as presidents of some well known universities. For example Tosun Terzio\u{g}lu, the head of the Turkish Mathematical Society between 1989 and 2008,  was a young colleague of Arf at METU, and  he knew Arf through Arf's friendship with his father Naz{\i}m Terzio\u{g}lu\footnote{He was a student of Caratheodory and  is best remembered for his leadership in the mathematical community and also finding lost works of  ibn al-Haytham. He was also among the team of mathematicians, including Cahit Arf, who in 1948 founded the Turkish Mathematical Society. He served as the head of the society between 1956 and 1976 but his administrative power was exercised mostly during his terms as dean and rector at several universities. After his death T\"{U}B{\.I}TAK honored his name with the Science Service Award in 1982.}. He did not refrain from taking responsibilities as an administrator. He served as the head of T\"{U}B{\.I}TAK and together with his team was instrumental in the smoothing out of obstacles that hindered the development of science in Turkey as well as opening up of some strategically important institutes. Terzio\u{g}lu then became the president of Sabanc{\i} University. Similarly Attila A\c{s}kar, another mathematician from the new generation, did not hesitate to take responsibility as the president of Ko\c{c} University. Both T. Terzio\u{g}lu and A. A\c{s}kar are recipients of T\"{U}B{\.I}TAK Science awards in mathematics. They had learned not only from Arf's preachings  and achievements but also equally well  from his  failures.

Arf's lifetime work was on the Riemann hypothesis where his efforts did not receive the same range of acceptance as his other work. However, it was his attitude towards mathematics through his relation with this problem that was instructional to his students. For example, from time to time rumors would arise telling   that the hypothesis was proven by some  mathematician somewhere. His response would always be cheerful. He would say that he is happy to finally learn the solution. He would make plans of clarifying those points which would inevitable require further explanation. A few days later he would be even more joyous because he had learned that the solution was using the same ideas that he was employing all along. Then as has always been the case, there would be a gap in the proof and Arf would again be happy that he still had some work to do. This altruistic approach to mathematics, in a \emph{publish or perish} world of the academia, with its sins and blessings affected most Turkish mathematicians.

After he resigned from  METU in 1980 he resided as a part time researcher at T\"{U}B{\.I}TAK's Gebze research center.
His inability to follow written mathematics  gradually caught up with him. In the 80's mathematics had long passed the threshold  where pure determination, perseverance and intelligence would solve problems. These traits were still necessary but no more sufficient. These on the other hand were Arf's strongest and may be even only assets in his mathematical endeavors. He was no longer  being fed by intellectually challenging seminars either. He started to slip back. Not to mention that old age was also pulling its tolls.  These however never stopped him from working on mathematics and chasing his ideas. He would work all day and complain later that he was having difficulties in recalling which indices he assigned to which variables.

When he could no longer travel to Gebze he would mostly stay at his home in Bebek, overlooking the Bosphorus, and occasionally visit Bo\u{g}azi\c{c}i University. He died on 27 December 1997 peacefully in his sleep. He was at that time survived by his wife, daughter and four grandchildren.
At the time of his death the President of Turkish Republic,
Mr. S\"{u}leyman Demirel, an engineer by education,
was a former student of Arf from the day's of Arf's residence at {\.I}stanbul Technical University. Arf was probably the only mathematician in the country
whose death was on prime time news.  His passion for scientific truth
was praised by all.

Posthumously, Arf was the subject of
three documentaries. One of these documentaries was prepared by private enterprise and was called ``Smorg", a mythological name referring to those who do all the work themselves. The other two documentaries were prepared by the state television TRT under the titles ``Cumhuriyete Kanat Gerenler (Those who kept the Republic under their Wings)" and ``I\c{s}{\i}kla Yaz{\i}lm{\i}\c{s} \"{O}yk\"{u}ler (Tales written with Light)".

Cahit Arf received numerous national awards and honors.
He received in 1948  the {\.I}n\"{o}n\"{u} Award,
in 1974  T\"{U}B{\.I}TAK Science Award and in 1988 Mustafa Parlar Foundation's Service in Science and Honor Award. Later in 1993  he was elected as a honorary member of T\"{U}BA, then the newly established Turkish Academy of Sciences.  He received honorary doctorates from {\.I}stanbul University and Karadeniz Technical University in 1980, and from METU in 1981. His only international recognition was being chosen as a member of Mainz Academy of Sciences and Literature.

Arf published his works in Turkish, English, French,  German and Italian. He was  also a fluent speaker in the first four of these languages.

A few years after his death METU Mathematics Department and the Mathematics Foundation initiated Arf Lectures\footnote{See http://www.matematikvakfi.org.tr/tr/arf-lectures}. Each year a prominent mathematician is invited to METU to give a talk at the mathematics auditorium which now bears Arf's name. The  series which started in 2001, and the speaker list includes as of 2013,  Gerhard Frey, Don Zagier, David Mumford, Robert Langlands, Peter Sarnak,
Jean-Pierre Serre, Hendrik Lenstra, Gunter Harder, Ben Joseph Green, John W. Morgan, Jonathan Pila and David E. Nadler.

In 2005 T\"{U}BA, Turkish Academy of Sciences, published a comprehensive biography of Arf, \cite{tosun}.

\section{Mathematical Activity of Arf}
We now proceed to record Cahit Arf's contributions to mathematics, mainly through those works which bear his name. For his role in the historical development of Turkish mathematics during the early years of the republic era after the collapse of the Ottomans  the reader may consult the accounts of Saban who was one of the collaborators of Arf, \cite{saban1,saban2}.
\subsection{Hasse-Arf Theorem}
In July of 1938 Arf completed his doctoral studies at G\"{o}ttingen University under the supervision of Helmut Hasse. His dissertation was submitted to Crelle's journal in November of that year and was published in the first issue of 1940. In this paper \cite{carf1} Arf generalized a work of Hasse \cite{hasse1,hasse2} and Artin \cite{artin}, and this generalization is commonly referred to as Hasse-Arf theorem. Nowadays the main reference for the modern version of Hasse-Arf theorem is  Serre's book {\it Local Fields}  , \cite[pages 76 and 93]{serre-local}. Among the numerous work done on Hasse-Arf theorem we find an interesting entry from 1972 by Ikeda, \cite{ikeda1}. Ikeda was a Turkish mathematician of Japanese origin and at that time he was in Middle East Technical University together with Arf. It must have been an exciting work for both gentlemen. In fact a brief but exact description of Hasse-Arf theorem  can be found in Ikeda's short article \cite{ikeda2} found in the appendix of Arf's collected works \cite{CW} published by the Turkish Mathematical Society to mark the occasion of Arf's 80th birthday.  The book however was prepared and published in 1988, two years ahead of time. Incidentally it was again Ikeda who wrote a detailed survey of Arf's work in number theory after Arf's death, \cite{ikeda3}.

Here we give a brief description of Hasse-Arf theorem following mainly the exposition given in Serre's book \cite{serre-local}. Let $L/K$ be a Galois extension with $G=G(L/K)$ its Galois group, where $K$ is a field complete under a discrete valuation $\nu_K$ and $L$ is a finite separable extension of $K$.
Let $\nu_L$ be the extension of the valuation $\nu_K$. Denote the valuation rings of $L$ and $K$ respectively by $A_L$ and $A_K$. Choose an element $x\in A_L$ generating $A_L$ as an $A_K$-algebra.
For any integer $i\geq -1$ define
\[ G_i=\{ s\in G\; |  \nu_L(s(x)-x)\geq i+1 \; \} . \]
These $G_i$'s define a decreasing filtration of normal subgroups of $G$ with $G_{-1}=G$ and $G_i=\{ 1\}$ for sufficiently large $i$. The group $G_i$ is called the $i$-th ramification group of $G$.

To define an upper numbering of the ramification groups we first define the function
\[ \phi(u)=\int_0^u \frac{dt}{(G_0:G_t)} \]
with the convention that $G_t$ is $G_i$ where $i$ is the smallest integer greater than or equal to $t$. The map $\phi$, which is continuous, piecewise linear and increasing, is a homeomorphism of the half plane $[-1,\infty)$ onto itself. The upper numbering is now defined as
\[ G^{\phi(u)}=G_u. \]
The lower numbering is convenient to study subgroups of the Galois group whereas the upper numbering is suited for quotients.

We say that a real number $v$ is a jump for the filtration if $G^v\neq G^{v+\epsilon}$ for all $\epsilon>0$. This is equivalent to saying that $\phi(u)$ is a jump whenever $G_u\neq G_{u+1}$.

Suppose that the group of quaternions $G$ is the Galois group of a totally ramified extension $L/K$. We then have a filtration where $G_0=G_1=G$, $G_2=G_3=C(G)$, the center of $G$, and $G_4=\{ 1\}$. In this case $v=3/4$ is a jump, see \cite[p77]{serre-local}.

Hasse-Arf theorem states that the jumps are always integers when $G$ is Abelian.

A more recent and a most beautiful account of Hasse-Arf theorem and its significance within the context of Artin-Hasse correspondence can be found in Roquette's\footnote{Peter J. Roquette, now a professor emeritus at Universit\"{a}t Heidelberg, was also a student of Hasse. His graduation date was 1951.} article \cite{roquette}.

There is also an interesting exposition tracing a chain of ideas from Langlands through Arf back to Hasse and his mathematical circle around this particular problem, see \cite{roquette2}

\subsection{Arf Invariant}
There are many different ways of `multiplying' two vectors to obtain a number and naturally mathematicians are interested in classifying all such possible multiplications. The classification naturally depends on the number system used as the target of the multiplication. Our understanding of this problem was almost complete while Arf was working on his dissertation. In 1938 when Arf successfully completed his graduate work he wanted to return to Turkey and serve idealistically at Kastamonu high school but Hasse convinced him to stay in G\"{o}ttingen for one more year. The problem Hasse suggested to Arf was the classification of quadratic forms over a field $K$ of characteristic 2.

Hasse's choice of this problem was very timely. A few years ago Witt\footnote{Witt and Arf  met and befriended each other. Witt later came to  {\.I}stanbul  to visit Arf. Witt family and Arf family had a close friendship besides the mathematical affinity of the fathers.} in \cite{witt} had given a complete system of invariants for quadratic forms except for those over a field of characteristic 2. Arf's mission was to fill this gap, which he did  in two consecutive papers, \cite{carf1,carf2}. In the former article he gives a full set of invariants of this quadratic form over a field of characteristic 2. The second article is a generalization of this result to formal power series rings in characteristic 2. Vanishing or non vanishing of a particular invariant defined in the first article is now known as the Arf invariant of that quadratic form.

Pioneer work, in most cases, is clumsy in mathematics. Sometimes a mathematician gives a long and elaborate proof of a notoriously difficult problem which has resisted the efforts of many prominent mathematicians over the ages, yet before long numerous short and transparent proofs start to appear in print. The pioneer has actually not solved the problem but has demolished a psychological barrier between mortals and the solution. The same phenomenon reoccurred with the definition of Arf invariants. Nowadays there is a {\it democratic} way of deciding on the Arf invariant $Arf(q)$ of a quadratic form $q$ on a vector space $V$ over $F_2$. Arf invariant is decided by the rule
\[ Arf(q)=\left\{ \begin{array}{ll}
1 & \mbox{if~} q(x)=1 \mbox{~for more than half the elements $x$ of $V$,} \\
0 & \mbox{otherwise.}
\end{array} \right. \]

Another way of defining the Arf invariant of the quadratic form $q$ on a $\Zz_2$-vector space $V$ is as follows:
\[ Arf(q)=\sum_{i=0}^{n}q(a_i)q(b_i)\mod 2 \]
where $a_i$'s and $b_j$'s form a symplectic basis of $V$ in the sense that $(a_i,b_j)=\delta_{ij}$, $(a_i,a_j)=(b_i,b_j)=0$ where $(x,y)=q(x+y)+q(x)+q(y)$ is the associated bilinear form of $q$ over $\Zz_2$. This is precisely the invariant Kervaire uses as an obstruction in constructing a triangulated manifold which admits no differentiable structure, see \cite{kervaire}.

Sophie Germain is attributed for the observation that ``Algebra is nothing but written geometry and geometry is nothing but algebra with figures." Following this adagio  it is expected that there should be a geometric interpretation of the Arf invariant. We look no further than the de~Rham interpretation of differential forms on differentiable manifolds. If the manifold has $2n$ real dimensions, then two differential forms of length $n$ can be multiplied to obtain a volume form which after integration over the manifold yields a number, and this is a perfect setting for the application of a quadratic form. Tensoring these operations with $F_2$ takes us into the realm of Arf invariants. This invariant has played an important role in some classification problems in algebraic and differential topology.

A detailed description of the algebraic aspects of Arf invariant can be found in \cite{ikeda2,ikeda3} and its topological significance  can be found in \cite{onder}.

Recently an error was discovered in the proof of a lemma in Arf's article \cite{carf2}. F. Lorenz and P. Roquette corrected this error and clarified the scope of Arf's theorem accordingly in two papers, see \cite{LR1,LR2}. A beautiful survey of the history and the humane elements surrounding this story is narrated in \cite{LR3}. The authors however note for the curious that all the facts and methods which are used in the new proofs were already present in Arf's original paper.

There is an interesting story behind the above story itself. Roquette asked in 2010 a question about finding a quadratic space with certain properties which would constitute a counterexample to Arf's theorem. Roquette was not a follower of the discussion list  {\tt math}{\it overflow} so his question was posted on the forum by a friend with the nickname \emph{KConrad}. This question was forwarded to Detlev Hoffman who was not a follower of the list either. His answer was posted back to the list by a user with the nick \emph{Skip}. A historical mathematical problem was thus settled by the mathematical community on the internet collaboratively, see \cite{mathoverflow}. As Roquette himself says ``Mathematics is done by people", and we may add ``irrelevant of whether they follow web forums or not." An entry about this whole story on Arf's theorem can also be found in Roquette's upcoming book \cite{roquette3}.

To further  emphasize the continuing effect of Arf invariant in current research, we remind here the recent work of M. A. Hill, M. J. Hopkins and D. C. Ravenel who solved the Kervaire-Arf invariant problem,  see \cite{hill1,hill2,hill3}.

\subsection{Arf Rings}
During one of  Patrick Du Val's seminars at {\.I}stanbul University he described his theory of the resolution of space curve branches. The resolution depended on some geometrically significant yet computationally mysterious invariants. He also published an article on this in {\.I}stanbul University journal, \cite{duval1}, the Turkish abstract of which can be easily recognized as being written by Arf in his own typical style. Following this seminar and the article Arf gave a full description of these invariants, using only the properties of the local ring of the branch at the singularity. The idea was to measure how far away these rings are from a {\it good} ring. This {\it good} ring was later called an Arf ring by Lipman and the process of finding the smallest Arf ring containing a ring was called the Arf closure of the ring, \cite{lipman}. This ring theoretic process
has a counterpart in the world of natural numbers. The orders of the elements in the completion of the local ring of the singularity form a numerical semigroup. The smallest semigroup which has the same Arf closure, similarly defined as in the ring case, as this semigroup has a minimal set of generators which are called the Arf invariants of the branch. The ideas used in this solution are heavily algebraic despite the geometric nature of the original problem. This had prompted Du Val to write a geometric description of Arf's article, \cite{duval2}.

Here we give a brief technical description of Arf rings and closures. A branch singularity in affine $n$-space is given formally by the parametrization
\beq
X_1&=&\phi_1(t) \\
&\vdots& \\
X_n &=&\phi_n(t)
\eeq
where each $\phi_i(t)$ is a formal power series in $t$ with no constant term. The smallest of the orders of the $\phi_i(t)$'s is called the multiplicity of the branch at the origin. Define
\[ H=k[[\phi_1(t),\dots,\phi_n(t)]], \]
the ring formally generated by the $\phi_i(t)$'s in the formal power series ring $k[[t]]$.

The blowing up of the branch at the origin produces the parametrization
\beq
X_1&=&\phi_1(t) \\
X_2&=&\frac{\phi_2(t)}{\phi_1(t)}-\lambda_1\\
&\vdots& \\
X_n &=&\frac{\phi_n(t)}{\phi_1(t)}-\lambda_n
\eeq
assuming that the order of $\phi_1(t)$ is the smallest. Here the constants $\lambda_j$'s are used to bring the singularity to the origin if it is not already there.
Define as above
\[ H_1=k[[\phi_1(t),\frac{\phi_2(t)}{\phi_1(t)}-\lambda_1, \dots, \frac{\phi_n(t)}{\phi_1(t)}-\lambda_n]]. \]

Here again we call the smallest order among the power series used to parameterize the branch the multiplicity at this stage. Continuing in this way we obtain a sequence of positive integers, called the multiplicity sequence.

The problem is to read off the multiplicity sequence from the local ring $H$.  If each $\phi_i(t)$ is of the form $t^{n_i}$ for some positive integer $n_i$, then the multiplicity sequence can be immediately derived from this initial data. However the obvious difficulty in the general case is that some of the trailing terms in some $\phi_i(t)$ may eventually become leading terms and there is no telling when this will happen.

For example if at some stage we have $\phi_1(t)=t^n$ and $\phi_2(t)=t^n+t^{n+m}$, then after the blow up we will have $\dis \frac{\phi_2(t)}{\phi_1(t)}-1=t^m$ and this $m$, which was never observed as the order of the series before, is now the leading term and may contribute to the multiplicity sequence.

To overcome this difficulty we must find a way of filling in the gaps.  For this we make the definition that a subring $H$ of $k[[t]]$ is an \emph{Arf ring} if for every positive integers $m$ and every $S_m\in H$ of order $m$, the set
\[ \{ \frac{S}{S_m}\in k[[t]] \; | \; \ord{S}\geq m \} \]
is closed under multiplication. This means simply that no surprises will emerge during the later stages of the blowing up process. The ring $k[[t]]$ is itself clearly an Arf ring. This then makes the following definition meaningful. The Arf closure $\arf H$ of the ring $H$ is the intersection of all the Arf rings containing $H$.

Similarly we define an Arf semigroup. A semigroup $G\subset\Nn$ is called an Arf semigroup if for all $m\in G$ the set
\[ \{ n-m\in \Nn \; | \; n\in G, \; n\geq m \} \]
is a semigroup. The semigroup $\Nn$ being Arf, the Arf closure $\arf G$ of $G$ is defined as the smallest Arf semigroup containing $G$.

Towards constructing $\arf H$ from $H$ we first observe that the ring $k+\phi_1(t)\cdot H_1$ contains $H$ and is contained in $\arf H$. Moreover it is clear after a second's thought that $\arf H=k+\phi_1(t)\cdot \arf H_1$. Repeat this process with $H_1$ replacing $\phi_1(t)$ with any element in $H_1$ of smallest nonzero order. This process eventually stops when $H_r$ is an Arf ring for some $r$.

Now we have a finite algorithm to recover the multiplicity sequence from the ring $H$. Starting with $H$ we first construct its Arf closure $\arf H$. Then we find the smallest semigroup $G$ whose Arf closure is equal to the semigroup of orders of $\arf H$. The generators of $G$ are called the Arf characters (or Du Val characters, see \cite{duval1}) of the branch. These are then used with a modified version of the Jacobian algorithm to give explicitly the multiplicity sequence. The Jacobian algorithm is the equivalent of Euclidian algorithm in finding the greatest common divisor of more than two integers. The use of this slowed-down version of the Jacobian algorithm and the geometric description of characters are first mentioned in \cite{duval1}.

There is a short description of Arf rings in an appendix to Arf's collected works, \cite{sertoz1}. Another but detailed survey of Du Val characters and Arf invariants for a singular curve branch was prepared for the 85th birthday of Arf, \cite{sertoz2}.

There are numerous applications of Arf rings and related invariants of a curve singularity in geometry. The applications range from complex geometry of the curves to algebraic codes. A relatively recent attempt to apply these ideas in higher dimensions was given by Campillo and Castellanos in \cite{campillo} where instead of using one dimensional objects they use one codimensional objects and  from this setting construct a one dimensional algebraic ring where Arf's ideas apply.

\subsection{Work on Elasticity}
His work on elasticity did not receive wide acknowledgement abroad but was highly praised and awarded in Turkey. His articles on this subject include \cite{carf6,carf9,carf10,carf12,carf13,carf16}. Despite the fact that these results did not attract attention then, the topic is currently active in particular about the structure of composite materials which are used in space science where any reduction in weight is of paramount importance.  There is  a brief note in the appendix to his complete works about Arf's work on elasticity, \cite{tezer}. A much more detailed survey of Arf's work on plane elasticity can be found in \cite{suhubi}.

\subsection{Lectures on Algebra}
This is Arf's only book, \cite{arfkitap}. Published in 1947 it is with today's standards a very ambitious piece of work.  Arf had  enormous energy for doing calculations and for that reason he was contend when he could reduce a problem to a stage where mere calculations would solve it. Apart from a few sections which  we would today skip in search of stronger theories to attack the problem, the book can be conveniently used  as an algebra course for undergraduate students. It has eight chapters. The first few chapters deal with rings, rings of polynomials in one and several variables. In particular meticulous detail is given to prime factorization of polynomials both in one and several variables. After treating group theory he enters into linear and algebraic equations and field extensions. There is a chapter on Galois theory where after treating standard material he talks about constructions with compass and ruler, and devotes two sections to describing the solutions of a cubic and quartic polynomial. The last chapter is on algebraic geometry where he starts with the resultants and proves Hilbert's Nullstellensatz. He then treats irreducibility of algebraic varieties in detail in terms of the defining system of polynomials and ends the book with a treatment of projective varieties and Bezout's theorem. It is a 450 page book full of exercises, 130 in total. The book reflects both the intellectual depth and the vigor of Arf's approach to mathematics.

In the introduction he mentions that the book started as lecture notes of algebra course he gave at {\.I}stanbul University between the years 1941-1944. He says that has paid special attention to the fact that a proof should make the reader feel why and how certain concepts are related besides being rigorous and correct. He claims that the book has no originality apart from the way he derives determinants from Cramer's rules and his treatment of the irreducibility of algebraic varieties. He has also published his method on determinants separately, see \cite{carf7}.

\section{Arf as an Intellectual Force}
In 1967 Langlands came to Ankara to spend a year at Middle East Technical University. He occupied an office next to Cahit Arf's. After a successful conference in {\.I}zmir, he took a walk with Arf along the streets of the city. During their conversation there Arf called Langlands a genius. In \cite{langlands} Langlands gives a beautiful narration of that incidence and says\footnote{Langlands wrote this in Turkish. The translation belongs to the author.}
\begin{quote}
Even though  I am not a genius, or may be only because of that, I never forgot that compliment. The only person in my life who ever called  me a genius was Cahit Arf. Even if you don't count the other reasons, this alone is enough for him to still hold a special place in my memories.
\end{quote}
In fact for exactly the same reason Cahit Arf still holds a special place in the memories of many Turkish mathematicians! Cahit Arf was  aware of his influence and would compliment a young mathematician at the most needed time and still manage to give the impression that the compliment is unique. This well monitored boost of morale worked miracles.

His support for working young scientists, be them mathematicians or not, was beyond mere verbal encouragements. At one time during his prime time in {\.I}stanbul University, he used his influence to have himself appointed as the head of the theoretical physics chair with the sole purpose of protecting Feza G\"{u}rsey\footnote{Feza G\"{u}rsey (1921-1992) was a prolific Turkish mathematician and physicist who worked in Turkey until 1974 and in Yale University after that time. He received numerous professional awards, an Einstein medal and an Oppenheimer  award being among them.}, a promising young member of the department, who at that time was having difficulties with the university administration in getting permission for extended foreign contact with his colleagues abroad.

His sincerity in his approach to science was a role model. He was not aiming at publishing no matter what comes. His foremost concern was to understand. He would demonstrate this in admitting his failures and errors openly. This set  ethical standards among his students which prevail today.

His standing up for his beliefs during the difficult times of METU and his doing this with self respect was another model. He demonstrated a living example that an intellectual can not and should not seek refuge in an ivory tower and be responsible only for science. He had sharp political views but knew how far he would push them without hurting other peoples' prides. He was hence respected even by those who strongly criticized his proposals.

His mathematical achievements and personal integrity overshadowed his shortcomings. For example he did not have students in the full sense of the word. He had some Ph.D. students but he did not create a school of thought in mathematics in Turkey. Towards the end of his life he realized this and made attempts to describe life long mathematical projects for the younger generations but they were not realistic and it was too late. He should have tried this when he was at his prime and still actively in touch with the mathematical world. But then he had always preferred to leave people alone to work on whatever comes close to their hearts. This approach with its sins and blessings marked the  Cahit Arf era in Turkey.

\section{A historical relic}

\begin{center}
\resizebox{15cm}{!}{\includegraphics{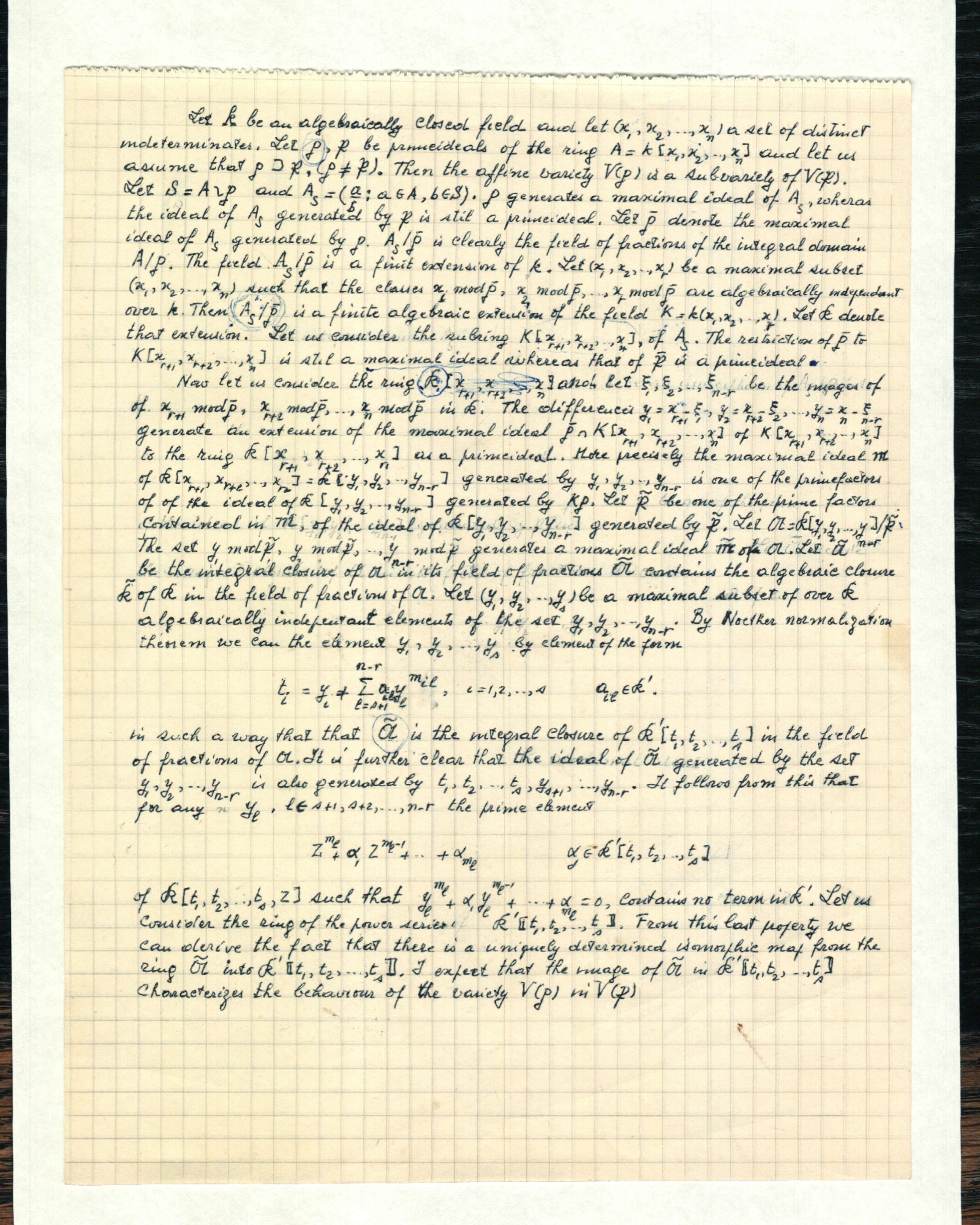}}
\end{center}

Between 1985 and 1988 I was a junior research assistant at TUBITAK's Gebze research institute. I was fresh out of graduate school and enjoyed company with Cahit Arf. We were mainly discussing the geometry of space curve branches but our discussions were restricted to discussing the geometric and historical concepts around the topic. We could not go into technical details for two reasons. Arf at that time was 75+ years old and he had not read any technical article during the past decade or two. But he had an enviable grasp for mathematical concepts. And he remembered very vividly his encounters and discussions with mathematicians long dead or just text book names for me. Mainly he would talk and explain the meanings, the actual meanings, of geometric ideas which are not openly written in any book.

Sometimes a young PhD student from abroad would come to the institute during summer times to meet Arf and talk about mathematics. Arf would invite him to give an impromptu informal talk about his  dissertation problem\footnote{One such student was Bet\"{u}l Tanbay, the head of Turkish Mathematical Society at the time of the writing of this article. However she says that she was duely warned beforehand that Arf would ask many questions, so she came well prepared and did not go through the shock described above.}. When the young student started to give his seminar Arf would interrupt ceaselessly to ask the meaning of every term used, even the most elementary ones. The young student giving the seminar would be utterly disappointed. I could read on his face the cry ``Is this the legendary Arf? He knows nothing!"  But after half an hour Arf would stop the seminar and would start to make comments about the proper way of attacking the problem  and discuss in detail its relation to other well known problems. His remarks were based only on the definitions he had heard during that seminar. Yet his remarks would startle the young student. His  face would blush with excitement. Often the student would cry ``That is exactly the heart of the technical problem". I witnessed this scenario several times during my stay at Gebze. This did not happen to me since the first week I was there Arf asked for a copy of my dissertation, read it at home and we discussed certain conceptual details later on. I would discover his secret later.

Sometimes after our discussions he would write, as usual on a small paper, with erasable ink pen, what he believed our topic of discussion needed in order to proceed properly. Here is one of those papers. I hardly remember why he wrote this and what we talked about it.  But judging from the scribblings  on the paper, I understand that we discussed on this paper. As usual the paper is not to be taken verbatim. It is the idea that matters, no matter how bad or wrong it is expressed.

Once I remember asking Arf about a certain statement in his paper on Arf rings, \cite{carf8}. He stared at the paper for a minute. Then he stared at me for a second. Smiled. And asked ``When was it that this was written?" It was written 40 years ago and I did not answer. I knew that I had to figure it out myself. Now I give a transcription of the paper below, line by line, and should the reader contemplate of asking me a question about its contents, I will simply remind the gentle reader that it was almost 30 years ago that I discussed this with Arf.

\begin{enumerate}
\item Let $\mathbf{k}$ be an algebraically closed field and let $(x_1,x_2,\dots,x_n)$ a set of distinct
\item indeterminates. Let $\pp$, $\qq$ be primeideals of the sing $A=k[x_1,x_2,\dots,x_n]$ and let us
\item assume that $\pp \supset \qq$, $(\pp\not=\qq)$. Then the affine variety $V(\pp)$ is a subvariety of $V(\qq)$.
\item
Let $S\setminus \pp$ and $A_S=\left( \frac{a}{b}:a\in A, b\in S\right)$. $\pp$ generates a maximum ideal of $A_S$, whereas
\item
the ideal of $A_S$ generated by $\qq$ is still a primeideal. Let $\bar{\pp}$ denote the maximal
\item
ideal of $A_S$ generated by $\pp$.  $A_S/\bar{\pp}$ is clearly the field of fractions of the integral domain
\item
$A/\pp$. The field $A_S/\bar{\pp}$ is a finite extension of $\mathbf{k}$. Let $(x_1,x_2,\dots,x_r)$ be a maximal subset [of]
\item
$(x_1,x_2,\dots,x_n)$ such that the classes $x_1 \mod \bar{\pp}$, $x_2 \mod \bar{\pp}$,\dots, $x_r \mod \bar{\pp}$ are algebraically independent
\item
over $\mathbf{k}$. Then $A_S/\bar{\pp}$ is a finite algebraic extension of the field $\mathbf{K}=\mathbf{k}(x_1,x_2,\dots,x_r)$. Let $\mathit{k}$ denote
\item
that extension. Let us consider the subring $\mathbf{K}[x_{r+1},x_{r+2},\dots,x_n]$ of $A_S$. The restriction of $\bar{\pp}$ to
\item
$\mathbf{K}[x_{r+1},x_{r+2},\dots,x_n]$ is still a maximal ideal whereas that of $\bar{\qq}$ is a primeideal.
\item
Now let us consider the ring $\mathit[x_{r+1},x_{r+2},\dots,x_n]$ and let $\xi_1,\xi_2,dots,\xi_{n-r}$ be the images of
\item
of $x_{r+1}\mod \bar{\pp},x_{r+2}\mod\bar{\pp},\dots,x_n\mod \bar{\pp}$ in $\mathit{k}$. The differences $y_1=x_{r+1}-\xi_1,y_2=x_{r+2}-\xi_2,\dots, y_n=x_n-\xi_{n-r}$
\item
generate an extension of the maximal ideal $\bar{\pp}\cap \mathbf{K}[x_{r+1},x_{r+2},\dots,x_n]$ of $\mathbf{K}[x_{r+1},x_{r+2},\dots,x_n]$
\item
to the ring $\mathit{k}[x_{r+1},x_{r+2},\dots,x_n]$ as a primeideal. More precisely the maximal ideal $\mathfrak{m}$
\item
of $\mathit{k}[x_{r+1},x_{r+2},\dots,x_n]=\mathit{k}[y_1,y_2,\dots,y_{n-r}]$ generated by $y_1,y_2,\dots,y_{n-r}$ is one of the primefactors
\item
of of the ideal of $\mathit{k}[y_1,y_2,\dots,y_{n-r}]$ generated by $\mathbf{K}\pp$. Let $\tilde{\qq}$ be one of the prime factors
\item
contained in $\mathfrak{m}$, of the ideal of $\mathit{k}[y_1,y_2,\dots,y_{n-r}]$ generated by $\tilde{\qq}$.
Let  $\al=\mathit{k}[y_1,y_2,\dots,y_{n-r}]/\tilde{\qq}$,
\item
The set $y_1 \mod \tilde{\qq}, y_2 \mod \tilde{\qq}, \dots y_{n-r}\mod \tilde{\qq}$ generates a maximal ideal $\tilde{\mathfrak{m}}$ of $\al$. Let $\tilde{\al}$
\item
be the integral closure of $\al$ in its field of fractions[.] $\tilde{\al}$ contains the algebraic closure
\item
$\tilde{\mathit{k}}$ of $\mathit{k}$ in the field of fractions of $\al$. Let $(y_1,y_2,\dots,y_s)$ be a maximal subset of over $\mathit{k}$
\item
algebraically independent elements of the set $y_1,y_2,\dots,y_{n-r}$. By Noether normalization
\item
theorem we can the elements $y_1,y_2,\dots, y_s$ by the elements of the form
\item
\[ t_i=y_i+\sum_{\ell=s+1}^{n-r}a_{i\ell}y_{\ell}^{m_{i\ell}}, \; i=1,2,\dots,s, \; a_{i\ell}\in \mathit{k}'.  \]
\item
In such a way that that $\tilde{\al}$ is the integral closure of $\mathit{k}'[t_1,t_2,\dots,t_s]$ in the filed
\item
of fractions of $\al$. It is further clear that the ideal of $\tilde{\al}$ generated by the set
\item
$y_1,y_2,\dots,y_{n-r}$ is also generated by $t_1,t_2,\dots,t_s,y_{s+1},\dots,y_{n-r}$. It follows from this that
\item
for any $y_{\ell}$, $\ell\in s=1,s+2,\dots,n-r$ the prime element
\item
\[ Z^{m_{\ell}}+\alpha_1Z^{m_{\ell}-1}+\cdots+\alpha_{m_{\ell}} \hspace{1cm} \alpha_j\in\mathit{k}'[t_1,t_2,\dots,t_s] \]
\item
of $\mathit{k}[t_1,t_2,\dots,t_s,Z]$ such that $y_{\ell}^{m_{\ell}}+\alpha_1y_{\ell}^{m_{\ell}-1}+\cdots+\alpha_{m_{\ell}}=0$, contains no terms in $\mathit{k}'$. Let us
\item
consider the ring of the power series $\mathit{k}'[[t_1,t_2,\dots,t_s]]$. From this last property we
\item
can derive the fact that there is a uniquely determined isomorphic map from the
\item

ring $\tilde{\al}$ into $\mathit{k}'[[t_1,t_2,\dots,t_s]]$. I expect that the image of $\tilde{\al}$ in $\mathit{k}'[[t_1,t_2,\dots,t_s]]$
\item
characterizes the behaviour of the variety $V(\pp)$ in $V(\qq)$
\end{enumerate}

\section{Concluding Reminiscences}
Cahit Arf was a major force on the face of scientific scene in Turkey. In a country who needed immediate development, he stood up for values which would pay off later.

His publication list is relatively short compared to today's standards, however the contents of his work and in particular their impact were impressive. Arf had the good fate of seeing his name  attached to his mathematical structures during his lifetime and this elevated him to a statue where it was easy for the younger generations to simply revere and follow him without questioning. He however used his influence to inflict free thought and independent thinking. As a result of this, the younger generations of mathematicians kept their Arf admiration separate from their independent evaluation of his achievements and failures. This is not a small achievement for a man who failed to form his school of thought in mathematics.

Though he was educated mainly abroad, his sense of belonging to his country played a major role in keeping many mathematicians working here in their own country. This, again with its cons and pros, is his mark on Turkish scientific scene.

Judging from the fact that he was much more productive when he was in the company of other active mathematicians, for example during Du Val's visit to {\.I}stanbul or Arf's own visit to Maryland, it is sometimes speculated that he would have produced enormously had he stayed at G\"{o}ttingen after his dissertation or at Maryland where in one year he had already produced two theorems. But then common sense reminds that in 1940's Turkey was a better place to be and Arf then made the right decision. The Maryland scenario presents another dilemma. Can we predict even for the sake of a thought experiment that a warm and amicable person such as Arf would perform better in an impersonal and harsh competitive environment? He was happy in the relaxed atmosphere of the Turkish scene of science then and did produce enough to have his name attached to mathematical concepts, and to more than one in that! At the end of such idle speculations one is tempted to conclude that one should produce where one is happy. He was happy in Turkey and did remarkably well.

Cahit Arf was very popular. In fact in many occasions we would see elementary school children herded into his office at Gebze by their teachers to show them what a real mathematician looks like. He would endure such visits with good humor and chat with the children in his usual witty style. He always had good stories to tell, even to the very young.

Needless to say, he was a legend in his own time and enjoyed it thoroughly.

\section{Some important dates in Arf's life}
\begin{tabular}{ll}
1910 & Born in Thessaloniki \\
1912 & Family moved to {\.I}stanbul \\
1926-1928 & St Louis Lyc\'ee years \\
1928-1932 & \'Ecople Normale Superior years \\
1932-1933 & Teacher at Galatasaray Lyc\'ee \\
1933 & Joined {\.I}stanbul University Mathematics Department \\
1937 & Arrived at G\"{o}ttingen to work with Helmut Hasse \\
1938 & Received his PhD at G\"{o}ttingen \\
1940 & Returned to {\.I}stanbul University \\
1943 & Promoted to professorship \\
1948 & Received {\.I}n\"{o}n\"{u} Science Award \\
1948 & Among the founders of Turkish Mathematical Society \\
1949-1950 & Visited The University of Maryland  \\
1955 & Promoted to ordinarius professor  \\
1956 & Membership at Mainz Academy \\
1962 & Retired from {\.I}stanbul University and joined Robert College \\
1963 & Among the founders of TUBITAK \\
1964 & Head of Science division at TUBITAK \\
1964-1966 & Visited the Institute of Advanced Studies at Princeton \\
1966-1967 & Visited University of California, Berkeley \\
1967-1980 & Worked at METU Mathematics Department \\
1974 & Received TUBITAK Science Award \\
1974 & Symposium in Arf's honor at Silivri \\
1980 & Received a honorary doctorate from {\.I}stanbul  University \\
1980 & Received a honorary doctorate from Karadeniz Technical University \\
1980 & Started to part time work at Gebze,  Fundamental Sciences Institute \\
1981 & Received a honorary doctorate from METU \\
1983-1989 & Head of Turkish Mathematical Society \\
1988 & Received Mustafa Parlar Science Award \\
1993 & Chosen as an honorary member of Turkish Academy of Sciences \\
1994 & Received Commandeur des Palmes Acad\'emiques from France \\
1997 & Died in {\.I}stanbul
\end{tabular}

\section{Ph.D. Students of Arf}
\begin{itemize}
\item Alt{\i}nta\c{s} G\"{u}le\c{c} (B\"{u}ke), 1949, Karakteristi\u{g}i 2 olan bir $K$ cismi \"{u}zerinde t\^{a}rif edilen, indeksi 3 olan kom\"{u}tatif-asosiyatif nilpotent cebirlerin invaryantlar{\i}, [Invariants of commutative associative nilpotent algebras of index 3 defined over a field of characteristic 2], ({\.I}stanbul University)\footnote{For a complete list of all degrees conferred at {\.I}stanbul University mathematics department between the years 1937 and 1981 see \cite{ozemre}.},
see \cite{buke}.
\\
\item Selma Soysal, 1949, Hilbert uzay{\i}nda birimin par\c{c}alan{\i}\c{s}{\i} ve baz{\i} s{\i}ralanm{\i}\c{s} projeksiyon operat\"{o}rleri c\"{u}mlelerine dair, [On the decomposition of unity on Hilbert spaces and on some ordered projection operator spaces], ({\.I}stanbul University),
see \cite{soysal}.
\\
\item Hil\^{a}l Pamir, 1953, Abstrakt kompleksler yard{\i}m{\i}yla topolojik uzay in\c{s}\^{a}s{\i} hakk{\i}nda, [On the construction of topological spaces using abstract complexes], ({\.I}stanbul University).
\\
\item Kemal \"{O}zden, 1954, Kabuklar{\i}n hesab{\i} hakk{\i}nda, [On the calculation of shells], ({\.I}stanbul University).
\\
\item Bediz Asral, 1955, Parabolik bir denklem i\c{c}in  Cauchy probleminin \c{c}\"{o}z\"{u}m\"{u} hakk{\i}nda, [On the solution of the Cauchy problem for parabolic equations], ({\.I}stanbul University), see \cite{asral}.
\\
\item \c{S}\"{u}k\^{u}fe Yamant\"{u}rk, 1955, {\.I}\c{c}i k\"{u}re \c{s}eklinde oyulmu\c{s} el\^{a}stopl\^{a}stik bir k\"{u}b\"{u}n emniyet katsay{\i}s{\i}n{\i}n t\^{a}yini, [Determination of the safety constant of a cube with a spherical cavity],
({\.I}stanbul University).
\\
\item Bahattin \c{S}aml{\i}, 1957, Dengede s\^{a}bit gerilmeli iki serbest s{\i}n{\i}rl{\i} d\"{u}zlem el\^{a}stik b\"{o}lgeler hakk{\i}nda, [On the elastic plane regions with two free boundary and constant tension at equilibrium],
({\.I}stanbul University).
\\
\item Ne\c{s}'et Ay{\i}rtman, 1961, Yal{\i}nkat fonksiyonlar{\i}n Taylor katsay{\i}lar{\i} hakk{\i}nda, [On the Taylor coefficients of univalent functions],
({\.I}stanbul University).
\\
\item \"{O}zdem \c{C}elik, 1961, Aritmetik lineer transformasyonlarla Dirichlet teoreminin katsay{\i}lar cismi sonlu olan bir rasyonel fonksiyonlar cismine te\c{s}mili, [On the generalization of Dirichlet theorem using arithmetical linear transformations to a rational function field whose constant field is finite],
({\.I}stanbul University), see \cite{celik}.
\\
\item G\"{u}ltekin B\"{u}y\"{u}kyenerel, 1972, Fibre bundles and normed rings, (METU).
{\it \small METU Library Call Number: QA611.B99}
\\
\item Halil {\.I}brahim Karaka\c{s}, 1974, On some invariants of algebraic curves,  (METU).
{\it \small METU Library Call Number: QA565.K18}
\\
\item Murat Kirezci, 1981, On the structure of universal non-IBN Rings $V_{n,m}$, (METU).
{\it \small METU Library Call Number: QA251.5.K58}
\\
\item Azize Hayfavi, 1983, On the L\"{o}wner Theory, (METU).
{\it \small METU Library Call Number: QA371.H419}
\end{itemize}

\section{Acknowledgements}
I thank my colleagues and friends Do\u{g}and \c{C}\"{o}mez , Metin G\"{u}rses, Mefharet Kocatepe, and  Erg\"{u}n Yal\c{c}{\i}n for reading  this manuscript and making encouraging remarks and correcting several errors. Tosun Terzio\u{g}lu told me  intimate anecdotes about Arf over a long phone talk, most of which found their way into the final version. Erdo\u{g}an \c{S}uhubi read an earlier version and made numerous suggestions which I implemented. Giocomo Saban kindly read an earlier version and  helped me correct a few historical facts. I am grateful to Nurettin Ergun and Yusuf Avc{\i} who helped me complete the list of Arf's students. I thank Peter Roquette who read the manuscript thoroughly and made numerous comments and suggestions which enriched and enhanced the final version. Finally, it was Reinhard Siegmund-Schultze's idea to dig up a handwritten note of Arf and share it on this manuscript.

It is also a pleasure for me to express special thanks to Kirsten Ward and Caroline Yal\c{c}{\i}n who read the previous versions and made meticulous remarks about language, style and presentation.

Despite the best efforts of my friends certain inaccuracies might have crept into the manuscript, for which I am  solely responsible.

\newpage

\renewcommand{\refname}{Scientific Publications of Arf}

\begin{center}
\resizebox{6cm}{!}{\includegraphics{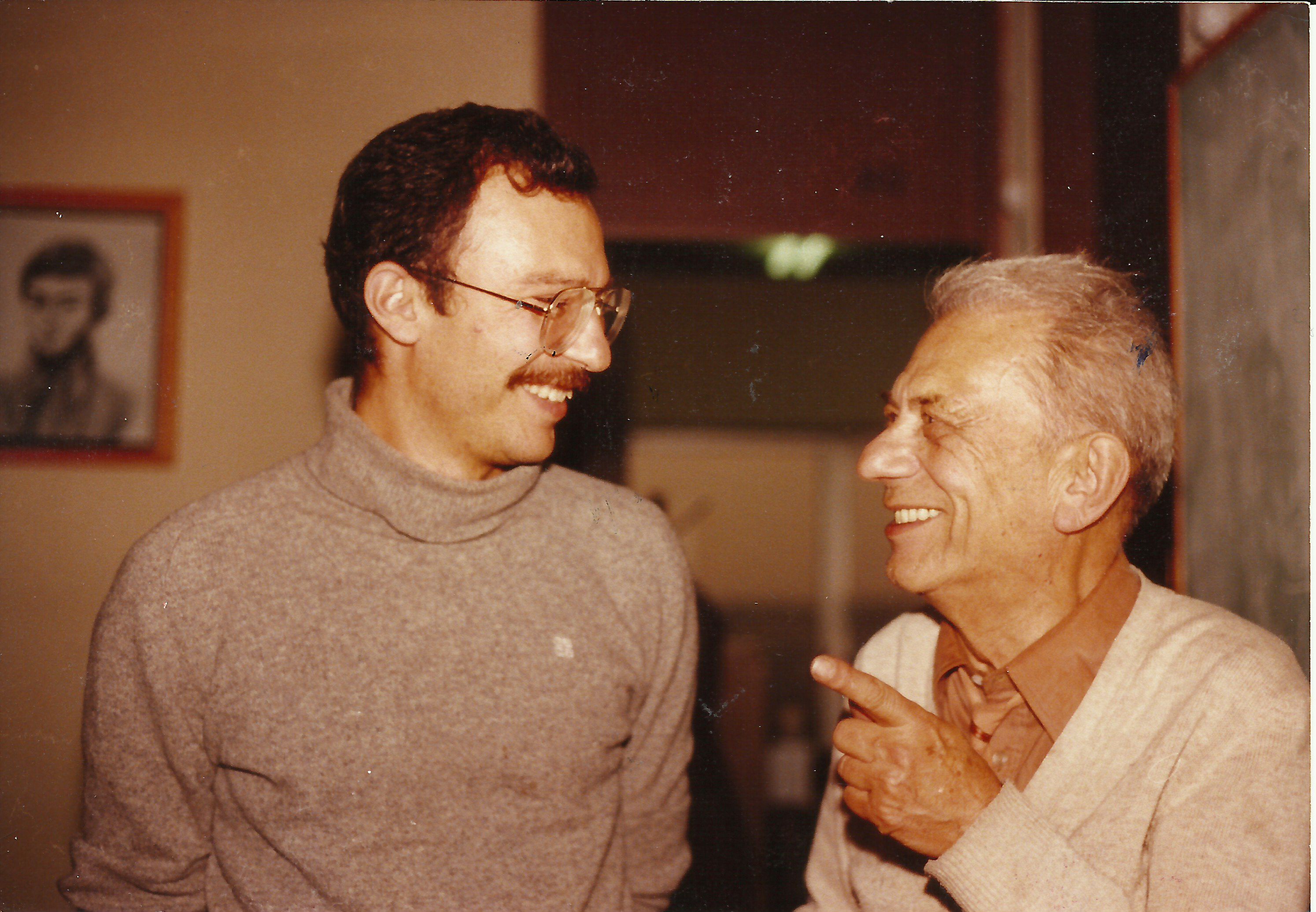}}
\resizebox{6cm}{!}{\includegraphics{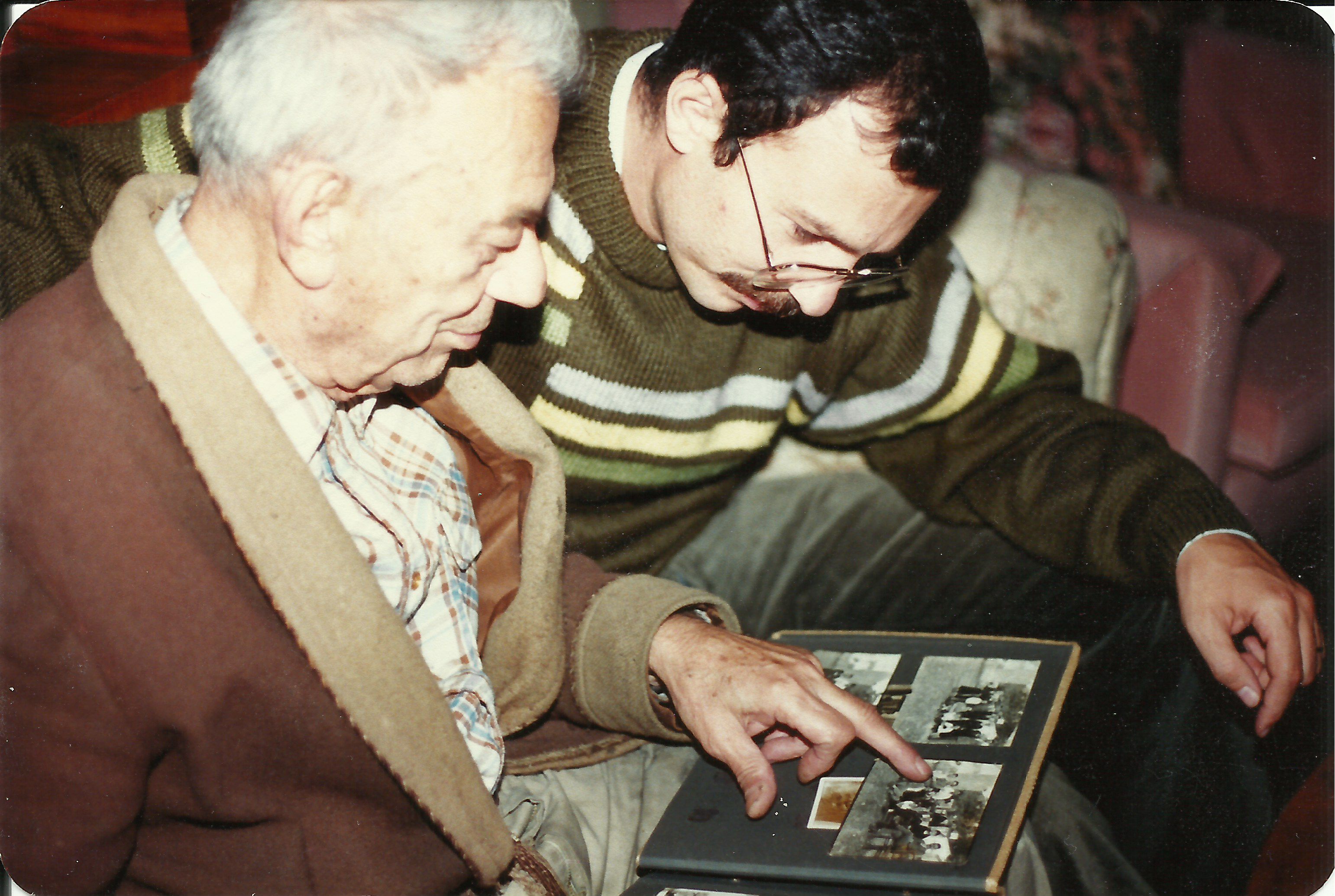}} \\
\makebox[6cm][l]{With Arf in 1985 in Gebze.} \makebox[6cm][l]{With Arf in 1988 in his Bebek home.}
\end{center}


\begin{thebibliography}{10}

\bibitem{carf1}
C.~Arf.
\newblock Untersuchungen \"uber reinverzweigte {E}rweiterungen diskret
  bewerteter perfekter {K}\"orper.
\newblock {\em J. Reine Angew. Math.}, 181:1--44, 1940.\\

\bibitem{carf2}
C.~Arf.
\newblock Untersuchungen \"uber quadratische {F}ormen in {K}\"orpern der
  {C}harakteristik 2. {I}.
\newblock {\em J. Reine Angew. Math.}, 183:148--167, 1941. \\

\bibitem{carf3}
C.~Arf.
\newblock Untersuchungen \"uber quadratische {F}ormen in {K}\"orper der
  {C}harakteristik 2. {I}{I}. \"{U}ber aritmetische {A}equivalenz quadratischer
  {F}ormen in {P}otenzreihenk\"orpern \"uber einem vollkommenen {K}\"orper der
  {C}harakteristik 2.
\newblock {\em Rev. Fac. Sci. Univ. Istanbul (A)}, 8:297--327, 1943. \\
\bibitem{carf4}
C.~Arf.
\newblock Un th\'eor\`eme de g\'eom\'etrie \'el\'ementaire.
\newblock {\em Rev. Fac. Sci. Univ. Istanbul (A)}, 12:153--160, 1947. \\

\bibitem{carf5}
C.~Arf.
\newblock {Sur les zones de Brillouin.}
\newblock {\em Rev. Fac. Sci. Univ. Istanbul, Ser. A}, 12:161--163, 1947.\\

\bibitem{carf6}
C.~Arf.
\newblock Sur la d\'etermination des \'etats d'\'equilibre d'un milieu
  \'elastique plan admettant des fronti\`eres libres \`a tensions constantes.
\newblock {\em Rev. Fac. Sci. Univ. Istanbul (A)}, 12:309--344, 1947.\\

\bibitem{arfkitap} C. Arf, {\em Cebir Dersleri} [Lectures on Algebra], Ankara University Faculty of Science Publications, No: 46 general, 9 mathematics, {\.I}brahim Horoz Publishing House,  442 pages, {\.I}stanbul, 1947,\\

\bibitem{carf7}
C.~Arf.
\newblock Sur la d\'efinition des d\'eterminants.
\newblock In {\em Universit\'e d'Istanbul. Facult\'e des Sciences. Recueil de
  m\'emoires comm\'emorant la pose de la premi\`ere pierre des Nouveaux
  Instituts de la Facult\'e des Sciences}, pages 9--20.
  Istanbul, 1948.\\

\bibitem{carf8}
C.~Arf.
\newblock Une interpr\'etation alg\'ebrique de la suite des ordres de
  multiplicit\'e d'une branche alg\'ebrique.
\newblock {\em Proc. London Math. Soc. (2)}, 50:256--287, 1948. \\


\bibitem{carf9}
C.~Arf.
\newblock Sur l'existence de la solution d'un probl\`eme d'\'elasticit\'e.
\newblock {\em Rev. Fac. Sci. Univ. Istanbul (A)}, 14:75--85, 1949. \\

\bibitem{carf10}
C.~Arf.
\newblock Sur l'identit\'e d'un probl\`eme de fronti\`ere libre en
  \'elasticit\'e avec un probl\`eme d'\'ecoulement.
\newblock {\em Bull. Tech. Univ. Istanbul}, 4(1):1--4 (1952), 1951. \\

\bibitem{carf11}
C.~Arf.
\newblock On the methods of {R}ayleigh-{R}itz-{W}einstein.
\newblock {\em Proc. Amer. Math. Soc.}, 3:223--232, 1952. \\

\bibitem{carf12}
C.~Arf.
\newblock On the determination of multiply connected domains of an elastic
  plane body, bounded by free boundaries with constant tangential stresses.
\newblock {\em Amer. J. Math.}, 74:797--820, 1952. \\


\bibitem{carf13}
C.~Arf.
\newblock Sur un probl\`eme de fronti\`ere libre d'\'elasticit\'e
  bidimensionelle.
\newblock {\em Bull. Tech. Univ. Istanbul}, 5:13--16 (1953), 1952. \\

\bibitem{carf14}
C.~Arf.
\newblock Remarques \`a propos d'un m\'emoire de {K}. {E}rim.
\newblock {\em Rev. Fac. Sci. Univ. Istanbul (A)}, 19:45--54, 1954. \\

\bibitem{carf15}
C.~Arf.
\newblock On a generalization of {G}reen's formula and its application to the
  {C}auchy problem for a hyperbolic equation.
\newblock In {\em Studies in mathematics and mechanics presented to Richard von
  Mises}, pages 69--78. Academic Press Inc., New York, 1954. \\

\bibitem{carf16}
C.~Arf.
\newblock Sur les fronti\`eres libres \`a tensions constantes d'un milieu
  \'elastique plan en \'equilibre.
\newblock {\em Rev. Fac Sci. Univ. Istanbul. S\'er. A.}, 19:119--132, 1954. \\

\bibitem{carf17}
C.~Arf.
\newblock Sur le th\'eor\`eme de {R}eiss.
\newblock {\em Rend. Mat. e Appl. (5)}, 14:181--191, 1954. \\

\bibitem{carf18}
C.~Arf.
\newblock Eine explizite {K}onstruktion der separablen {H}\"ulle eines
  {P}otenzreihenk\"orpers.
\newblock {\em Abh. Math. Sem. Univ. Hamburg}, 19:117--126, 1955. \\

\bibitem{carf19}
C.~Arf.
\newblock \"{U}ber den {S}atz von {D}ubourdieu.
\newblock {\em Abh. Math. Sem. Univ. Hamburg}, 20:112--114, 1955. \\

\bibitem{carf20}
C.~Arf.
\newblock \"{U}ber ein {A}nalogon des {R}iemann-{R}ochschen {S}atzes in
  {Z}ahlk\"orpern.
\newblock {\em Akad. Wiss. Mainz. Abh. Math.-Nat. Kl.}, 1957:293--328, 1957. \\

\bibitem{carf21}
C.~Arf, K.~{\.I}mre, and E.~{\"O}zizmir.
\newblock On the algebraic structure of the cluster expansion in statistical
  mechanics.
\newblock {\em J. Mathematical Phys.}, 6:1179--1188, 1965. \\

\bibitem{carf22}
C.~Arf and G.~Saban.
\newblock Due dimostrazioni elementari del teorema di {C}artan e {D}ieudonn\'e.
\newblock {\em Riv. Mat. Univ. Parma (2)}, 6:169--223, 1965. \\

\bibitem{carf23}
C.~Arf.
\newblock Sur la structure du groupe de {G}alois de la fermeture alg\'ebrique
  d'un corps de s\'eries de puissances sur un corps fini et les conducteurs
  d'{A}rtin.
\newblock In {\em Les Tendances G\'eom. en Alg\`ebre et Th\'eorie des Nombres},
  pages 27--35. \'Editions du Centre National de la Recherche Scientifique,
  Paris, 1966. \\

\bibitem{carf24}
C.~Arf.
\newblock The advantage of geometric concepts in mathematics.
\newblock In {\em Rings and geometry (Istanbul, 1984)}, pages 553--556. Reidel,
  Dordrecht, 1985.  \\

{\quad}

\textsc{\large References}

{\quad}

\bibitem{CA} \emph{Cahit Arf}, a booklet prepared by Middle East Technical University on the occasion of their conferring a honorary doctorate to Arf in 1981. \\

\bibitem{CW} C. Arf, {\it Complete Papers of Cahit Arf}, Turkish Mathematical Society, 1988. \\

\bibitem{arf-yok} C. Arf, \"{U}niversiteler yasas{\i} nedir ve ne olmal{\i}d{\i}r, [What is and what should be the (new) laws on university (education system)], \emph{Cahit Arf} booklet published by the Middle East Technical University on the occasion of his receiving a Honorary Doctorate from the university,  47-53, 1981. \\

\bibitem{arf-hatay} C. Arf, Iktisadi olaylar{\i}n incelenmesinde matematik metodlar{\i}n tatbiki hakk{\i}nda, [On the applications of mathematical methods on investigating economical affairs], \emph{Tenth University Week, Hatay}, {\.I}stanbul University Publications no: 716, Becid Publishing House, {\.I}stanbul, 1957.\\

\bibitem{artin} E. Artin, Gruppentheoretische Struktur der Diskriminante algebraischer Zahlk\"{o}rper, J. reine angew. Math., 164:1-14, 1931. \\

\bibitem{asral} Asral, B.,
On the solution of the Cauchy problem for parabolic equations,
Rev. Fac. Sci. Univ. {\.I}stanbul, S\'{e}r. A 21, 63-85, 1956, \\


\bibitem{buke}
B\"{u}ke, A., Untersuchungen \"{u}ber kommutativ-assosiativ und nilpotenten Algenren vom Index 3 und von der Charakteristik 2, Rev. Fac. Sci. Univ. {\.I}stanbul, S\'{e}r. A 19 Supplement, 1-145, 1954.
\\


\bibitem{campillo} A. Campillo\ and\ J. Castellanos, Valuative Arf characteristic of singularities, Michigan Math. J., 49:435-450, 2001. \\

\bibitem{celik}
\c{C}elik, \"{O}., Eine Anwendung der arithmetisch-linearen Transformationen auf den Dirichletschen Satz im algebraischen Funktionenk\"{o}rper, Rev. Fac. Sci. Univ. {\.I}stanbul, S\'{e}r. A 33, 77-110, 1968. \\


\bibitem{duval1} P. Du Val, The Jacobian algorithm and the multiplicity sequence of an algebraic branch,  Rev. Fac. Sci. Univ. {\.I}stanbul, 7:107-112, 1942. \\

\bibitem{duval2} P. Du Val, Note on Cahit Arf's ``Une interpr\'{e}tation alg\'{e}briques de la suite des ordres de  multiplicit\'e d'une branche alg\'ebrique", Proc. London Math. Soc., 50:288-294, 1949.\\

\bibitem{alp} Eden, A., Irzik, G, German mathematicians in exile in Turkey: Richard von Mises, William Proger, Hilda Geringer, and their impact on Turkish mathematics, Historia Mathematica, 39, 432-459, 2012. \\

\bibitem{gunning} R. C. Gunning and H. Rossi,
{\it Analytic Functions of Several Complex Variables}, Prentice-Hall Series in Modern Analysis, 1965.\\


\bibitem{hasse1} H. Hasse, F\"{u}hrer, Diskriminante und Verzweigungsk\"{o}rper relativ Abelscher Zahlk\"{o}rper, J. reine angew. Math., 162:169-184, 1930. \\


\bibitem{hasse2} H. Hasse, Normenresttheorie galoisscher Zahlk\"{o}rper mit Anwendungen auf F\"{u}hrer und Diskriminante abelscher Zahlk\"{o}rper, J. Fac. Sci. Tokyo, 2:477-498, 1934. \\



\bibitem{hill1}  Hill, Michael A. ;  Hopkins, Michael J. ;  Ravenel, Douglas C.  A solution to the Arf-Kervaire invariant problem.
 Proceedings of the G\"{o}kova Geometry-Topology Conference 2010,
 21--63, Int. Press, Somerville, MA,  2011.\\
		
\bibitem{hill2}  Hill, Michael A. ;  Hopkins, Michael J. ;  Ravenel, Douglas C.  The Arf-Kervaire problem in algebraic topology: sketch of the
 proof.
 Current developments in mathematics, 2010,
 1--43, Int. Press, Somerville, MA,  2011.\\
		
\bibitem{hill3}  Hill, Michael A. ;  Hopkins, Michael J. ;  Ravenel, Douglas C.  The Arf-Kervaire invariant problem in algebraic topology:
 introduction.
 Current developments in mathematics, 2009,
 23--57, Int. Press, Somerville, MA,  2010.\\

\bibitem{hironaka} H. Hironaka, Resolution of singularities of an algebraic variety over a field of characteristic zero, Annals of Math., 79 (1964), 109-326.\\

\bibitem{ikeda1} M. G. Ikeda, A generalization of Hasse-Arf theorem, J. reine angew. Math., 252:183-186, 1972. \\

\bibitem{ikeda2} M. G. Ikeda, Topics in algebraic number theory, Appendix in {\it The Collected Papers of Cahit Arf}, 410-412, 1988. \\


\bibitem{ikeda3} M. G. Ikeda, Cahit Arf's contribution to number theory and related fields, Turkish J. Math., 22:1-14, 1998. \\

\bibitem{kervaire} M. Kervaire, A manifold which does not admit any differentiable structure, Comment. Math., 34, 257-270, 1960.\\

\bibitem{langlands} R. Langlands, Benim bildi\u{g}im Cahit Arf (Cahit Arf that I know), {\tt http://www.sunsite.ubc.ca/DigitalMathArchive/Langlands/ \\ miscellaneous.html{\#}arf} \\

\bibitem{lipman} J. Lipman, Stable ideals and Arf rings, Amer. J. Math., 93:649-685, 1971. \\


\bibitem{LR1} Lorenz, F., Roquette, P., On the Arf invariant in historical perspective part 2, Mathematische Semesterberichte, 59, 2011. \\


\bibitem{LR2} Lorenz, F., Roquette, P., On the Arf invariant in historical perspective, Mathematische Semesterberichte, 57:73-102, 2010. \\

\bibitem{LR3} Lorenz, F., Roquette, P., Cahit Arf and his invariant, Mitteilungen der Mathematischen Gesellschaft in Hamburg", vol. XXX, 87-126, 2011. \\
    see http://www.rzuser.uni-heidelberg.de/~ci3/arf3-withpicture.pdf \\

\bibitem{mathoverflow} mathoverflow question 21720, Wanted: Quadratic Space in Characteristic 2 as a Counterexample to a Theorem of Arf, \\ see: http://mathoverflow.net/questions/21720/  \\


\bibitem{onder} T. \"{O}nder, Arf invariant and its applications in topology, Appendix in {\it The Collected Papers of Cahit Arf}, 413-415, 1988. \\

\bibitem{onder2} T. \"{O}nder, Bir Portre, Cahit Arf, \\
\verb+http://www3.iam.metu.edu.tr/matematikvakfi/cahitarf-onder.pdf\+
\\

\bibitem{ozemre} \"{O}zemre, A. Y., (Editor), {\.I}stanbul \"{U}niversitesi Fen Fak\"{u}ltesi'nde \c{C}e\c{s}itli Fen Bilimleri Dallar{\i}n{\i}n Cumhuriyet D\"{o}nemindeki Geli\c{s}mesi ve Milletleraras{\i} Bilime Katk{\i}s{\i}, [ The evolutions of various branches of sciences  during the Republic period at the Science Faculty of {\.I}stanbul University and their contributions to international science], {\.I}stanbul University Publications No: 3042, Science Faculty No: 176, {\.I}stanbul, 1982.\\

\bibitem{roquette} Roquette, P., On the history of Artin's L-functions and conductors, Seven letters from Artin to Hasse in the year 1930, Mitteilungen der Mathematischen Gesellschaft in Hamburg, Band  19$^\ast$ (2000) 5-50. \\
A revised version is avaliable at the author's site: \\
http://www.rzuser.uni-heidelberg.de/$\sim$ci3/lfunktio.pdf\\

\bibitem{roquette2} Roquette, P., Introduction of Langlands at the Arf Lecture, 11 November 2004, \\
    see: http://www.rzuser.uni-heidelberg.de/~ci3/chair.pdf \\

\bibitem{roquette3} Roquette, P., Contributions to the History of Number Theory in the 20th Century, European Mathematical Society, 2003. \\

\bibitem{saban1} G. Saban, Sviluppo Storico della Matematica nell'Impero Ottomano e durante i primi anni della Repubblica Turca, Bolletino della Unione Matematica Italiana, (8), 5-A:73-96, 2002. \\

\bibitem{saban2} G. Saban, Sviluppo Storico della Matematica in Turchia dalla Riforma dell'Universit\'{a} al 1997, Bolletino della Unione Matematica Italiana, (8), 5-A:257-292, 2002. \\

\bibitem{serre-local} J-P. Serre, {\it Local Fields}, Springer-Verlag Graduate Texts in Mathematics no 67, 1979. \\

\bibitem{sertoz1} S. Sert\"{o}z, On Arf Rings,
Appendix in {\it The Collected Papers of Cahit Arf}, 416-419, 1988. \\


\bibitem{sertoz2} S. Sert\"{o}z, Arf Rings and Characters,
Note di Matematica 14:251-261, 1994.\\

\bibitem{soysal}
Soysal, S., Sur une repr\'{e}sentation fonctionnelle d'une d\'{e}composition g\'{e}n\'{e}rale de l'unit\'{e} dans l'espace de Hilbert, Bull. Tech. Univ. {\.I}stanbul 2, No: 2, 65-86, 1949.\\


\bibitem{suhubi} E. \c{S}uhubi, Cahit Arf'{\i}n elastisite teorisine katk{\i}lar{\i}, [Cahit Arf's contribition to the theory of elasticity], \emph{Cahit Arf} booklet published by the Middle East Technical University on the occasion of his receiving a Honorary Doctorate from the university, 40-46, 1981. \\

\bibitem{tezer} C. Tezer, Arf's work in applied mathematics,
Appendix in {\it The Collected Papers of Cahit Arf}, 420-422, 1988. \\

\bibitem{tosun} T. Terzio\u{g}lu, A. Y{\i}lmaz, {\it Cahit Arf, `Anlamak' Tutkunu Bir Matematik\c{c}i}, [Cahit Arf, A Mathematician devoted to `understanding'], Turkish Academy of Sciences, Biographies Series no: 4, Ankara, 2005. \\

\bibitem{witt} E. Witt, Theorie der quadrischen Formen in beliebigen K\"{o}rpern, J. reine angew. Math., 176:31-44, 1937. \\



\end{thebibliography}
\end{document}